\documentclass[reqno,12pt]{article}

\usepackage[all,line,poly,rotate,ps,dvips]{xy}
\usepackage{paralist,comment}
\usepackage[pdftex]{graphicx} 
\SelectTips{cm}{}

\usepackage[italian]{babel}
\usepackage[T1]{fontenc}

\begin{document}

\title
{Francesco Severi: il suo pensiero matematico e politico prima e dopo la Grande Guerra}

\author{
{\normalsize\sc Ciro Ciliberto}
\vspace{-0.75mm}\\
{\small Dipartimento di Matematica}\vspace{-1.4mm} \\
{\small Universit\`a di Roma ``Tor Vergata''}\vspace{-1.4mm}\\
{\small Via della Ricerca Scientifica}\vspace{-1.4mm}\\
{\small 00133 Roma, Italy}\vspace{-1.4mm} \\
{\small e-mail: {cilibert@mat.uniroma2.it} }\\
{\normalsize\sc Emma Sallent Del Colombo}
\vspace{-0.75mm}\\
{\small Departament d'ELL-DCEM}\vspace{-1.4mm} \\
{\small Universitat de Barcelona}\vspace{-1.4mm}\\
{\small Pg. de la Vall d'Hebron, 171}\vspace{-1.4mm}\\
{\small 08035 Barcelona, Spain}\vspace{-1.4mm} \\
{\small e-mail: {emma.sallent@ub.edu}} 
}

\date{}

\maketitle

\begin{flushright}
\emph {Dedicato ad Umberto Bottazini \\ 
in occasione del suo 70esimo compleanno}
\end{flushright}

\tableofcontents

\section{Introduzione}

Francesco Severi (1879--1961) è una delle figure più rilevanti dell'ambiente matematico italiano e internazionale nella prima metà del XX secolo. Si tratta di una figura molto complessa con una presenza importante nella vita politica e culturale italiana.

La  produzione scientifica di Severi comprende pi\`u di 400 pubblicazioni di ricerca matematica, numerosi trattati e scritti vari. I suoi contributi, multiformi e di notevole levatura, spaziano in vari settori della geometria algebrica: variet\`a abeliane e quasi abeliane, studio e  classificazione delle superficie e variet\`a algebriche a meno di trasformazioni birazionali, geometria proiettiva delle variet\`a algebriche, geometria numerativa; notevoli sono anche alcuni articoli di teoria delle funzioni analitiche di pi\`u variabili complesse.  Vanno ricordati anche i suoi libri di testo, in cui cerca di coniugare l'aspetto intuitivo con un notevole rigore metodologico. 

Tra i suoi principali allievi ricordiamo (in ordine di tempo): Ruggiero Torelli (1884--1915), Giacomo Albanese (1890--1948), Annibale Comessatti, Giuseppe Gherardelli (1894--1944), Fabio Conforto (1909--1954), Beniamino Segre (1903--1977), Enzo Martinelli (1911--1999), Federico Gaeta (1923--2007).  

Numerosi sono stati gli interventi storici che lo riguardano, molti dei quali verranno citati in seguito. In questo articolo ci siamo proposti di accomunare i vari aspetti della personalità di Severi, quello matematico, quello politico, quello istituzionale accademico e quello filosofico, cercando di capire in particolare se gli eventi relativi alla Grande Guerra hanno costituito un punto di svolta nella vita e nella produzione scientifica di Severi.

Sembra emergerne una figura poliedrica e di grande complessità per la quale il periodo bellico e quello immediatamente successivo, coincidono con una sostanziale affermazione della sua personalità che lo portano ad assurgere a una funzione di guida politica della matematica italiana e a svolgere un ruolo importante anche nella politica culturale del paese.

Nel contempo, la sua produzione scientifica, rispetto al periodo anteguerra, risente di un'eccessiva fiducia autoreferenziale che lo porta a una minore cura del rigore e quindi a dei risultati talvolta errati.

Abbiamo suddiviso quindi il lavoro in due blocchi, il primo corrispondente al periodo antecedente la Grande Guerra e il secondo a quello successivo. Ciascun blocco prende in considerazione i vari aspetti che riguardano la personalità e l'opera di Severi.

\section{Prima della Grande Guerra}

Francesco Severi nacque ad Arezzo il 13 aprile 1879. La sua infanzia fu segnata dalla tragica morte del padre Cosimo, notaio, avvenuta il 4 gennaio 1888. Essa ebbe gravi ripercussioni economiche sulla numerosa famiglia (madre e cinque figli, altri quattro erano premorti).
Severi ricorda nella breve biografia che \`e premessa a \emph{Dalla Scienza alla fede} del 1959 (cfr. \cite {Sev59}), che inizi\`o a dare lezioni private all'et\`a di dodici anni, per aiutare la famiglia.
  
Dopo essersi iscritto al corso di laurea di ingegneria dell'Universit\`a di Torino, capì, principalmente per l'influenza dell'insegnamento di Corrado Segre (1863--1924), che la sua vera inclinazione era per la matematica pura e cambi\`o indirizzo di studi in tal senso. Si laure\`o nel 1900 con una tesi  su un argomento di \emph{geometria numerativa} \cite {Sev00}, tema su cui torn\`o ripetutamente nel seguito.

Dopo la laurea, sempre a Torino, diventò assistente di Enrico D'Ovidio (1843--1933) e tenne, da giovanissimo libero docente, il corso di ``Geometria Proiettiva e Descrittiva''. Successivamente si trasferì all'Universit\`a di Bologna (1902--1903) come assistente di Federigo Enriques (1871--1946)  e quindi all'Universit\`a di Pisa come assistente di Eugenio Bertini (1846--1933).  Nel 1904 ottenne la cattedra di ``Geometria Proiettiva e Descrittiva'' all'Universit\`a di Parma. Dopo un solo anno passò all'Universit\`a di Padova, dove diventò direttore della Scuola di Ingegneria. 

Fin da subito prolifico autore di importanti lavori di geometria algebrica, non tardò ad ottenere importanti riconoscimenti in ambito nazionale ed internazionale. Nel 1906 ottenne la Medaglia d'oro dll'Accademia dei XL e, l'anno successivo, insieme a Federigo Enriques, il Prix Bordin dell'Accademia delle Scienze di Parigi per i contributi alla classificazione delle \emph{superficie iperellittiche}, infine la Medaglia Guccia nel 1908 al Congresso Internazionale dei Matematici di Roma, la cui commissione giudicatrice si componeva di
Max Noether (1844--1921), Henri Poincar\'e (1854--1912) e Corrado Segre. Divenne socio dell'Accademia Nazionale dei Lincei nel 1910 e dell'Accademia delle Scienze di Torino nel 1918. 

Con l'entrata dell'Italia nella prima guerra mondiale Severi, socialista ma al contempo interventista convinto, si arruola volontario nell'artiglieria, dove ha modo di distinguersi: si occupa di \emph{fonotelemetria}\footnote{Determinazione della posizione di una sorgente sonora in base agli istanti di arrivo delle onde sonore in tre punti del terreno non allineati; usata in guerra per localizzare la posizione di artiglierie nemiche, mediante il \emph{fonotelemetro}: strumento costituito da sensori infissi nel terreno che, analizzando le vibrazioni sonore prodotte dai colpi di artiglieria, ne determina la posizione di origine.} e, al pari di vari altri matematici coinvolti nel conflitto\footnote{Ad esempio Guido Fubini (1879--1943),  Mauro Picone (1885--1977) , Antonio Signorini (1888--1963), Alessandro Terracini (1889--1968), Vito Volterra (1860--1940), cfr. \cite {Tanzi77, Tanzi81}. \`E proprio dall'impegno di Picone e dello stesso Volterra nei calcoli numerici necessari per le correzioni al tiro che nasce la prima idea per l'istituzione di un Istituto di Calcolo, che poi diverr\`a l'Istituto Nazionale per le Applicazioni del Calcolo istituito nel 1932 presso il CNR.}
di correzioni al {tiro di artiglieria} dipendenti dalle variazioni di temperatura e pressione dovute all'altitudine,  scrivendo anche un lungo articolo al riguardo, pubblicato alla fine del conflitto (cfr. \cite {Sev19, Tanzi81}).

Col suo assistente Annibale Comessatti (1886--1945), ufficiale di artiglieria decorato con due croci di guerra, Severi viene dislocato in Val Lagarina (Rovereto) e raggiunge il grado di capitano.

\subsection{Argomenti matematici}

Il periodo precedente la  prima guerra mondiale \`e probabilmente quello in  cui Severi d\`a i contributi pi\`u rilevanti e compiuti alla geometria algebrica. Dal 1900, anno della laurea di Severi a Torino, relatore Corrado Segre, al 1918, Severi pubblica 89 tra lavori scientifici, testi e interventi. 

I suoi contributi spaziano su diversi argomenti di geometria algebrica come la geometria delle varietà proiettive, la geometrica numerativa, la classificazione delle superficie algebriche. Diamo qui di seguito un rapido cenno degli argomenti trattati e dei risultati principali.

\subsubsection*{Geometria delle variet\`a proiettive} Di particolare interesse l'articolo \cite {Sev01} in cui si  introducono  e classificano le superficie oggi dette OADP (one apparent double point) e si d\`a una fondamentale caratterizzazione della superficie di Veronese di grado 4 nello spazio proiettivo di dimensione 5 aprendo un solco che ha portato Fyodor Zak negli anni 1980 a classificare le cosiddette \emph{variet\`a di Severi} (cfr. il libro \cite {Zak93}). 
Nell'articolo \cite {Sev06} viene dimostrata la prima istanza di quel che poi \`e diventato il \emph{teorema della sezione iperpiana} di Salomon Lefschetz (1884--1972), che afferma che ogni sezione iperpiana liscia di una variet\`a intersezione completa liscia $X$ di dimensione almeno 3 in uno spazio proiettivo \emph{eredita} varie propriet\`a topologiche notevoli di $X$ (cfr. \cite {Lef24}).  Nella memoria \cite {Sev15}  vengono studiate le \emph{variet\`a grassmanniane} che parametrizzano i sottospazi di data dimensione di uno spazio proiettivo fissato, dimostrandone alcune propriet\`a fondamentali, in particolare che il loro gruppo di Picard \`e generato dalla sezione iperpiana nella immersione di Pl\"ucker, ossia che ogni loro sottovariet\`a di codimensione 1 \`e intersezione completa con una ipersuperficie dello spazio ambiente.

\subsubsection*{Geometria numerativa} La geometria numerativa si occupa di \emph{contare} il numero di oggetti algebrico--geometrici che variano in una data famiglia e che verifichino date propriet\`a che ne determinino appunto un numero finito. Uno degli archetipi di questo genere di questioni \`e il \emph{problema delle cinque coniche} di Jacob Steiner (1796--1863): le coniche in un piano proiettivo dipendono da 5 parametri, e il problema consiste nel determinare il numero (3264) delle coniche tangenti a 5 coniche sufficientemente generali fissate.   

Fondamentale \`e la  memoria  \cite {Sev02}, in cui Severi si muove nella direzione della soluzione del quindicesimo dei 23 problemi posti da David Hilbert (1862--1943) nel testo della sua conferenza tenuta l'8/8/1900 al Congresso Internazionale dei Matematici a Parigi. Severi cerca infatti di porre le basi per una trattazione rigorosa di problemi numerativi classici e della fondazione del calcolo numerativo proposto da Hermann Schubert (1848--1911). Severi abborda il problema con l'introduzione di caratteri numerativi delle variet\`a che, in termini moderni, sono  \emph{classi di Chern} di opportuni \emph{fibrati vettoriali} (cfr. \cite {Ful84} per una trattazione moderna). 

 \subsubsection*{Classificazione delle superficie algebriche} La classificazione delle superficie era stata affrontata con grandi successi (che Severi fece sempre fatica ad accettare, attribuendosi talvolta il merito anche di taluni di questi) da Guido Castelnuovo e da Federigo Enriques fin dagli inizi degli anni '90 del XIX secolo, quindi circa dieci anni prima che Severi iniziasse ad occuparsene. Quando Severi entra in scena, la classificazione \`e gi\`a delineata, ma restano importanti problemi aperti. Tra i maggiori, vi sono  la classificazione delle superficie regolari con sistema canonico banale, che noi oggi chiamiamo $K3$ secondo una denominazione suggerita da Andr\'e Weil (1906--1998) (cfr. \cite [p. 546] {Wei80}), e di quelle che chiamiamo \emph{biellittiche} mentre classicamente venivano chiamate \emph{iperellittiche}, cio\`e superficie irregolari con sistema canonico non banale, ma un cui multiplo \`e banale. 
 
Severi affronta entrambi i problemi. La classificazione (non esaustiva) delle superficie iperellittiche, effettuata in collaborazione con Federigo Enriques, frutta ai due autori il Prix Bordin del 1907 (cfr. \cite {EnrSev08,  EnrSev09}). La vicenda della competizione al premio tra Enriques e Severi da una parte e  Giuseppe Bagnera (1865--1927) e Michele De Franchis (1875--1946) dall'altra \`e ben nota (cfr. \cite {Cil04}) e non \`e necessario trattenerci qui su questo punto (cfr. \cite {CilSer91} per approfondimenti). Lo studio delle superficie $K3$, del cui \emph{spazio dei moduli} Severi d\`a una prima descrizione, si trova nella memoria \cite {Sev08}. 

Uno degli argomenti pi\`u delicati nella teoria delle superficie algebriche \`e (ancora oggi) lo studio delle superficie \emph{irregolari}. \`E questo uno degli argomenti cui Severi dedic\`o molta attenzione in tutta la sua carriera.  Le sue ricerche prendono spunto dall'articolo \cite {Sev03}, che tratta un argomento suggeritogli da Enriques (cfr. \cite {Cil04}) e si sviluppano nel corso dei decenni successivi. La questione principale, mai risolta nell'ambito algebrico--geometrico dalla scuola italiana, \`e quella di mostrare la coincidenza di varie possibili definizioni (aritmetica, geometrica, analitica, topologica) del concetto di \emph{irregolarit\`a} di una superficie, che a sua volta \`e una possibile estensione alle superficie del concetto di \emph{genere} delle curve. La questione divenne nota come il \emph{problema della completezza della serie caratteristica}.  Anche qui la storia \`e ben nota e non \`e il caso di intrattenervicisi ancora, mentre rimandiamo a  \cite [\S 2.2] {BCP04} per una visione generale del problema. 

In un paio di importanti memorie \cite {Sev07, Sev09} del 1907--09 Severi si lancia anche su problemi inerenti lo studio di variet\`a di dimensione superiore a due,  discutendo della estensione a tali variet\`a dei  generi delle superficie, delle irregolarit\`a e  del Teorema di Riemann--Roch. 

\subsubsection*{Il teorema della base}  Si tratta qui del problema fondamentale dello studio del gruppo delle curve algebriche sopra una superficie modulo l'\emph{equivalenza algebrica} (oggi detto \emph{gruppo di N\'eron--Severi} della superficie), al quale Severi d\`a fin dalle prime sue ricerche un contributo essenziale e fortemente originale. Anche qui le sue ricerche prendono il via dall'articolo \cite {Sev03/2} (correlato a \cite {Sev03}), che di nuovo verte su un argomento suggeritogli da Enriques:  in esso Severi effettua uno studio aritmetico delle corrispondenze tra curve, ossia delle curve appartenenti alla superficie prodotto di due curve. Il teorema principale, che oggi si enuncia dicendo che il {gruppo di N\'eron--Severi} di una superficie \`e un gruppo abeliano finitamente generato, si trova in \cite {Sev06/2}.

Severi fu autore anche di diversi trattati. Tra questi si distinguono le \emph{Lezioni di geometria algebrica} \cite {Sev08/2}, apparse poi in traduzione tedesca in \cite {Sev21} con una serie di interessantissime appendici. Questa traduzione, come segnalato in \cite {Sev15/2}, era gi\`a pronta nel 1915, ma non fu pubblicata per lo scoppio della guerra. Di \cite {Sev15/2} fu accelerata la pubblicazione, anche su stimolo di Georg Zeuthen (1839--1920),  per il timore che non fosse pi\`u possibile farlo pi\`u tardi causa la guerra. Tra le appendici, specie le F) e G), basate in buona sostanza su \cite {Sev15/2},  hanno giocato un ruolo fondamentale nello sviluppo della teoria delle curve algebriche nel XX secolo, introducendo idee e problemi su alcuni dei quali non \`e stata ancora detta l'ultima parola (come, ad esempio, lo studio delle \emph{variet\`a di Severi} di curve singolari su superficie, lo studio della struttura birazionale dello spazio dei moduli delle curve, o quello delle famiglie di curve (iper)spaziali). Purtroppo molte delle asserzioni ivi contenute non sono corrette: queste Appendici segnano infatti la linea di confine tra il periodo anteguerra (caratterizzato dalla produzione di risultati importanti dimostrati con una buona cura per il rigore), e quello dopoguerra (caratterizzato pi\`u da congetture che da teoremi dimostrati in modo rigoroso).

\subsection{La politica, i rapporti accademici, la filosofia}

\subsubsection*{La politica}

Goodstein a Babbitt ben riassumono in \cite {GooBab11} l'itinerario politico di Severi fino alla prima guerra mondiale: 

\begin{quote}

[...] as a boy he took a keen interest in politics, following in particular the socialist movement, then on the rise in Italy. After being appointed professor of mathematics at Padua in 1905, Severi allied himself with the left--wing \emph{blocco popolare patavino}, which rewarded his allegiance by appointing him president of the municipal gas and water company. In 1910 he officially joined the Socialist Party and was quickly elected councilor for the commune of Padua and became the Socialist alderman for education. When World War I broke out, Severi sided with those who urged intervention on the side of Britain and France. He severed connections with the Socialist Party, which supported Italian neutrality, and quickly volunteered for military service once Italy entered on the Allied side in 1915.
\end{quote}

Per delucidare meglio il carattere e le  scelte di Severi, \`e opportuna qualche testimonianza pi\`u diretta e qualche approfondimento. \`E ad esempio interessante menzionare un articolo sulla ``Cronaca di Padova'' del 8/2/1914,
che testimonia del forte coinvolgimento di Severi col Partito Socialista e con le sue posizioni anti--interventiste. Si parla di una conferenza tenuta la sera prima a Padova da Alfredo Rocco (1875--1935)\footnote{Famoso giurista, al cui nome \`e legato il codice penale da lui varato e tuttora in vigore, nonch\'e quello di procedura penale, rimasto in vigore dal 1930 fino al 1988. Rocco fu professore ordinario all'Universit\`a di Parma (1906--1909), a quella di di Palermo (1909--1910), e alla Universit\`a di Padova (1910--1925), per poi, con un percorso simile a quello di Severi, trasferirsi definitivamente all'Universit\`a di Roma, ateneo di cui fu rettore dal 1932 al 1935.
Rocco, prima vicino al partito radicale, diventa nazionalista nel 1913 e pubblica nel gennaio 1914 l'opuscolo \cite {Roc14}, che contiene il programma teorico del movimento, dai toni imperialistici e antidemocratici che anticipano le posizioni del fascismo. Rocco fu un acceso sostenitore dell'ingresso dell'Italia nella Grande Guerra. Dal 1925 al 1932 fu Ministro di grazia e giustizia e affari di culto e promosse la codificazione penale del fascismo, firmando, come gi\`a detto,  il codice penale e quello di procedura penale.}, allora anch'egli professore all'Universit\`a di Padova, sul tema ``Il nazionalismo e i partiti politici''.  Rocco fa affermazioni del tipo:

\begin{quote}
[...] il nazionalismo tende al predominio della razza italiana, considera gli interessi della razza come preminenti ed assoluti. [...] fra le due forme, emigrazione e conquista, quella che pi\`u si addice ad un popolo forte \`e quest'ultima.
\end{quote}

Rocco viene allora interrotto dal grido di ``Abbasso la guerra!'' levatosi da un gruppo di socialisti  capitanati da Severi, il quale pretende di parlare in risposta a Rocco, ma viene fermato dai sostenitori di quest'ultimo con lo slogan di ``Domandi il permesso a Benedetto Croce!'', del cui significato derisorio capiremo il significato tra breve. 

Severi cambiò però poi idee politiche. Fu molto probabilmente a causa di serie divergenze con la direzione del Partito Socialista sull'intervento in guerra, che Severi usc\`i dal partito. Infatti \`e del 9/3/1915, un anno dopo l'episodio sopra ricordato, un suo lungo articolo, sul giornale ``L'Adriatico'' che riportiamo per intero in appendice (cfr. \S\ref{App}), in cui difende appassionatamente la sua nuova posizione interventista. 

L'articolo sull'``Adriatico'' ci pare interessante perch\'e delucida alcuni tratti salienti del carattere e dello stile di Severi, e, in ultima analisi, del suo atteggiamento ogni volta che si \`e trovato a cambiare bandiera.  Nella lettera si legge l'autodefinizione di ``uomo di studio pi\`u che uomo di parte'' volto a contribuire ``sia pure in modestissima misura, a ricondurre a quel rispetto reciproco delle idee e a quella concordia degli animi, di cui [...] c'\`e tanto bisogno''.  Rispetto reciproco che non pareva preoccuparlo all'epoca della conferenza di Rocco. Questo atteggiamento, di persona al di sopra delle parti, pur se attivo nella vita politica del paese, ritorner\`a nell'opuscolo \cite {Sev59} per giustificare la presa di distanza dal fascismo e l'adesione alla fede cattolica nel secondo dopoguerra. Ed \`e pure tipica dell'atteggiamento di Severi la giustificazione delle sue inversioni di rotta umane e politiche con ragionamenti apparentemente logici, ma che in sostanza sono volti a piegare le ideologie prima abbracciate alle sue nuove scelte. 
 
\subsubsection*{I rapporti accademici}
Dal punto di vista accademico, Severi assurge ben presto ad una posizione di prestigio, sebbene  dimostri, nel periodo prebellico, grande rispetto verso Castelnuovo ed Enriques, che egli sembra considerare suoi mentori. Ad esempio, in una lettera del 31/10/1904\footnote{Le sette lettere di Severi a Castelnuovo, che vanno dal 1903 al 1915, si trovano nel fondo Castelnuovo  dell'Accademia Nazionale dei Lincei curato da Paola Gario sul sito:

\noindent http://operedigitali.lincei.it/Castelnuovo/Lettere\_E\_Quaderni/menu.htm

A quel che ci risulta esse non sono state ancora pubblicate n\`e commentate, pur contenendo elementi di sicuro interesse. In esse salta agli occhi la grande deferenza (e all'inizio vera e propria devozione) di Severi che d\`a sempre del Lei a Castelnuovo, anche una volta diventato ordinario. I rapporti tra Castelnuovo e Severi rimasero, se non cordiali, almeno buoni, anche dopo il contrasto di Severi con Enriques, le vicende degli anni del Fascismo, leggi razziali incluse, e l'epurazione (Severi fu epurato proprio da Castelnuovo dall'Accademia dei Lincei, e riammesso dallo stesso poco dopo). Anche dal punto di vista scientifico i rapporti erano presenti, come testimonia l'interessantissima lettera di Castelnuovo a Severi del 26/11/1947 anch'essa nel citato fondo.  Castelnuovo ebbe a dire, in occasione del giubileo scientifico di Severi nel 1950, che Severi era ``uno dei maggiori matematici che l'Italia abbia prodotto negli ultimi cento anni'', cfr.  \cite {Cas51}. I rapporti Castelnuovo--Severi meriterebbero un esame approfondito.}, Severi ringrazia Castelnuovo per gli auguri per la vittoria nel concorso di Parma e dice:
\begin{quote} Accetto gli auguri per l'avvenire, come incitamento ad operare per mantenere ed accrescere la stima di coloro che mi hanno guidato finora con aiuti e con amichevoli consigli.\end{quote}

I rapporti amichevoli e le collaborazioni del primo periodo tra Enriques e Severi, si muteranno  in un aperto contrasto.  Severi viene menzionato nella corrispondenza Enriques--Castelnuovo \cite {BCG96} per la prima volta il 5/12/1901, a proposito di una sua possibile collocazione quale assistente di Enriques a Bologna, e vi compare con una certa regolarit\`a (anche se non spesso) fino al 24/8/1914, nell'ultima lettera della corrispondenza a noi giunta.  In questo periodo non vi \`e traccia di contrasti tra Enriques e Severi nella corrispondenza. Anzi ne emerge una  comunanza di interessi scientifici e accademici tra i due.  In particolare si nota la cura e  l'apprensione di Enriques per la carriera di Severi, riguardo anche al percorso alquanto travagliato del concorso a cattedra di cui Severi risult\`o vincitore  nel 1904, che lo port\`o a divenire ordinario a Parma.   \`E tuttavia interessante notare qualche lettera di Enriques riguardo ai primi esordi di Severi il cui tono non \`e esattamente elogiativo. Ad esempio si noti il tono ironico, tipico del maestro impaziente di fronte ad un allievo che egli giudica neghittoso, nella lettera del 3/6/1902:

 \begin{quote}
 Sapevo invece del caso Severi\footnote{Per l'opposizione di Cremona e di Dini, vi era il pericolo che a Severi non venisse concessa la libera docenza nel 1902, in quanto non aveva alcun lavoro di geometria descrittiva.}, ed anzi ho sulla coscienza di avergli suggerito tre o quattro argomenti elementari di Geometria descrittiva, qualcuno dei quali fu da lui scartato perch\'e non abbastanza elementare. Si trattava infatti, nientemeno, che della rigata cubica, e di distinguere sulla rappresentazione i due casi in cui \`e bilatera o unilatera!\footnote {Problema che poi Severi risolse in una nota del 1903, cfr. \cite {Sev03/3}, attribuendosene grande merito, cfr. \cite [Vol. I, p. 118] {Sev71-89}.} 
\end{quote}

\noindent Successivamente nella lettera del 31/10/1902  si riferisce che Bertini ha chiesto, per motivi di salute, a Enriques di cedergli come assistente Severi, in cambio di Arturo Maroni (1878--1966). Enriques dice:

\begin{quote}
A malincuore non ho potuto rifiutare il Severi in queste circostanze\footnote{Effettivamente Severi, dopo aver trascorso  l'anno accademico 1902--03 come assistente di Enriques a Bologna, passer\`a a Pisa quale assistente di Bertini per un anno, prima di vincere il concorso a cattedra ed essere chiamato a Parma.}. Per\`o alla $2^ a$ parte della domanda (prendermi Maroni) non mi sento di consentire.

Non mi par giusto di preferire ai nostri \emph{mediocri}, un altro mediocre, non privo forse di studio, ma mal dotato per l'insegnamento stante il suo carattere eccessivamente freddo. 
\end{quote} 

Negli anni 1907--1909 si concentra la collaborazione scientifica tra Enriques e Severi sulla classificazione delle superficie iperellittiche.  Ma le interazioni tra i due riguardo allo studio delle superficie irregolari inizia anche prima, con i suggerimenti di Enriques a Severi di studiare la superficie delle coppie di punti di una curva di genere positivo, e poi successivamente con la dimostrazione, che si riveler\`a sbagliata, da parte di Enriques, del teorema di completezza della serie caratteristica (cfr. \cite {Enr04}). Quest'ultimo risultato fu utilizzato indipendentemente da Castelnuovo \cite {Cas05} e da Severi \cite {Sev05} per dimostrare che i vari concetti di irregolarit\`a che possono definirsi su una superficie sono tutti equivalenti. Come evidenziato, ad esempio, in \cite {Bri04, Cil04}, le sintonie tra Enriques e Severi in questo primo periodo non si fermano ai comuni interessi scientifici, ma si allargano sul piano didattico, su quello accademico, con l'interesse dimostrato da Enriques per le sorti concorsuali di Severi (cfr. \cite [p. 39]{Cil04} e \cite [p. 548, lettera del 9/9/1904]{BCG96}) e su quello filosofico su cui torneremo tra breve. Tuttavia, come suggerito in \cite [p. 47] {Cil04}, \`e probabilmente proprio nel corso delle collaborazioni di questi anni che le differenze caratteriali e del modo stesso di concepire la ricerca matematica  -- sistematico e metodico per Severi, fortemente intuitivo per Enriques -- che i semi dei futuri contrasti  gettano le radici.

In ogni caso, una prima avvisaglia di tali contrasti va registrata gi\`a prima della grande guerra. Ne abbiamo notizia dalla corrispondenza di Severi con Vito Volterra, che nell'epoca di cui parliamo era il vero nume tutelare della matematica in Italia, conservata nel Fondo Volterra presso l'Accademia Nazionale dei Lincei\footnote{Si tratta di una ventina di lettere che vanno dal 1900 alla fine degli anni '20. La corrispondenza presenta un  carattere di grande rispetto da parte di Severi verso Volterra, cui viene sempre dato del Lei, perfino dopo il trasferimento di Severi a Roma. Anche queste lettere non sono state a nostro avviso pubblicate. Esse non trattano questioni matematiche, ma ve ne sono varie che hanno interesse di natura istituzionale.}.

Da un paio di lettere della primavera del 1909 emerge un dissapore  tra Severi, allora presidente della Mathesis\footnote{Associazione fondata nel 1895  da Rodolfo Bettazzi (1861--1941),  volta ai problemi dell'insegnamento della matematica, voluta e creata da un gruppo di insegnanti di scuola secondaria, cui ben presto si affiancarono molti illustri docenti universitari, alcuni dei quali, tra cui gli stessi Castelnuovo, Enriques e Severi, ne assursero poi alla guida.}, da un lato e Castelnuovo, Enriques e Giovanni Vailati (1863--1909) dall'altro. Il contrasto era dovuto al fatto che i tre erano stati nominati dal Comitato Centrale per l'Insegnamento della Matematica, formato da Felix Klein (1849 --1925), Alfred  Greenhill (1847--1927) e Henri Fehr (1870--1954),  creato durante il Congresso Internazionale dei matematici di Roma del 1908, a costituire una Delegazione col compito di indagare sulle condizioni dell'insegnamento della matematica in Italia. Severi contestava la nomina, in quanto la Mathesis, il cui compito istituzionale era appunto di occuparsi di questioni riguardanti l'insegnamento in Italia, non era stata consultata nel formare tale delegazione e i ``triunviri'' intendevano mantenere una posizione del tutto indipendente rispetto alla  stessa Mathesis. \`E interessante un passo della lettera del 14/4/1909:
\begin{quote}
... da parte mia desidero che questa questione resti nel campo obbiettivo e che non guasti la mia amicizia con Enriques e Castelnuovo. Da questo lato mi sembra che Castelnuovo consideri la cosa al suo giusto valore; mentre Enriques vuol dimostrarsi evidentemente irato e imbronciato. Me ne duole moltissimo, ma quel che devo fare, lo faccio lo stesso!
\end{quote}  

Riguardo ai rapporti con Volterra, \`e anche interessante ricordare, perch\'e ha a che fare con gli eventi bellici,  il caso delle celebrazioni per il settantesimo compleanno di Max Noether (1844--1921) occorso nel settembre 1914, cui si riferisce uno scambio epistolare dell'epoca.  Volterra fa sapere a Severi, che si occupa della cosa sul versante italiano (il che dimostra la posizione di prestigio da lui gi\`a assunta sul piano nazionale e internazionale), che intende partecipare alle onoranze a Noether. Severi gli spiega che, per mancanza di tempo, si \`e potuta organizzare soltanto una pergamena firmata dai soli geometri algebrici e non da tutti i matematici italiani (peraltro una pergamena in tal caso sarebbe stato troppo poco); comunque Severi lascia Volterra libero di decidere se apporre la firma o no.   Volterra decide di firmare, ma il 17/9/1914 la pergamena si trova ancora a Padova, Severi ha avuto difficolt\`a a spedirla causa la guerra e allora chiede a Volterra di aiutare a farlo attraverso un canale diplomatico (non si sa se ci\`o sia poi riuscito). La vicenda \`e interessante: si noti infatti che entrambi, Volterra e Severi, erano  interventisti e Volterra diverr\`a nel dopoguerra aspramente antitedesco, al punto da approvare l'allontanamento degli studiosi tedeschi da convegni internazionali.

\subsubsection*{Severi e la filosofia}

Nel 1906 Severi scrive la breve recensione \cite {Sev06/3} ai \emph{Problemi della Scienza} di  Enriques \cite {Enr06}, il testo 
che raccoglie la visione critico--positivista dell'autore. L'introduzione al libro inizia cos\`i:  \medskip

\begin{quote}
Una riflessione, {maturatasi durante il decennio fra il 1890 e il 1900}, ci ha condotto alla critica di taluni problemi che si riferiscono allo sviluppo logico e psicologico delle conoscenze scientifiche; i quali vengono qui trattati come ``problemi della Scienza''. 
\end{quote}\medskip

I punti cruciali  sollevati in \cite {Enr06} sono quelli tipici della filosofia Enriquesiana, e cio\`e  l'unit\`a dei saperi,  la critica dei fondamenti delle varie discipline scientifiche (la matematica \`e solo una di queste), la spiegazione psico--fisiologica di tali fondamenti, l'affermazione della validit\`a della geometria come scienza naturale, che nasce dall'esperienza sensoriale e fa parte integrante della fisica. 

Severi, come gi\`a osservato ad esempio in \cite [\S 2] {Bri04}, ben interpreta il punto di vista di Enriques, iniziando la sua  recensione con le parole:

\begin{quote} Il volume qui recensito [...] si raccomanda per la visione unitaria che subordina la specificit\`a dei saperi ad una comune prospettiva gnoseologica.
\end{quote}

Negli stessi anni Giovanni Gentile (1875--1944), insieme a Benedetto Croce (1866--1952), cominciava a vedere in Enriques filosofo un pericoloso incursore nel loro campo di competenza, specie in quanto propugnatore di una filosofia scientifica a loro del tutto estranea. Il dissidio, come si sa, esplose in modo aperto in occasione del IV Congresso Internazionale di Filosofia organizzato da Enriques a Bologna nel 1911, anche questa una storia gi\`a troppe volte scritta per ritornarci qui in dettaglio, anche perch\'e non del tutto pertinente al nostro tema. 

Nella sua lunga e articolata recensione \cite {Gent08} ai ``Problemi della Scienza'' su ``La Critica'', rivista diretta da  Croce, Gentile scrive che  si tratta di  ``vagheggiamenti di una filosofia scientifica''. Ed egli accomuna nella sua critica anche la ``Rivista di Scienza'' (che dal 1910 cambier\`a il nome in ``Scientia''), da poco fondata da Enriques e da Eugenio Rignano (1870--1930) con un programma ben chiaro di rinnovamento filosofico enunciato nel ``Programma'' che apre, alle pagine 1--2, il primo numero della rivista:
\begin{quote}

Contro codesti criterii ristretti intende reagire soprattutto il movimento nuovo di pensiero verso la sintesi; una Filosofia libera da legami diretti coi sistemi tradizionali, sorge appunto a promuovere la coordinazione del lavoro, la critica dei metodi e delle teorie, e ad affermare un apprezzamento pi\`u largo dei {problemi della Scienza}. Pel quale il particolarismo stesso viene compreso in un aspetto pi\`u adeguato nella interezza del processo scientifico.

\end{quote}

Gentile liquida tutto ci\`o come qualcosa che ``non pu\`o se non incoraggiare il dilettantismo scientifico''. 
Allo stesso tempo, la scarsa dimestichezza dei filosofi idealisti con i concetti espressi da Enriques \`e testimoniata dalle seguenti parole che si trovano in una lettera di Gentile a Croce datata pochi giorni prima della suddetta recensione (cfr. \cite{CG3}):

\begin{quote}
Spero di scrivere domani la recensione dell'Enriques, che \`e un libro che non so per che verso pigliare, per non dirne troppo male con {la paura di non aver capito}, per colpa mia, ci\`o che ci pu\`o essere di buono.
\end{quote}\medskip

Una successiva incursione di Severi nel campo filosofico, probabilmente nella speranza di emulare anche in ciò Enriques (cfr. \cite {Cil04}) \`e  del 1910, quando scrive l'importante articolo \cite {Sev10}. Qui egli, a dispetto delle critiche idealiste, sposa in pieno il positivismo critico di Enriques:
 
 \begin{quote}
 ... la geometria sa di quel che parla: del mondo fisico. Essa differisce dalla fisica soltanto nel metodo, prevalentemente sperimentale per l'una, deduttivo nell'altra. E anche il metodo perde il suo carattere deduttivo quando si tratta di scoprire. Ai confini della scienza, nelle posizioni di avanguardia, la severit\`a logica \`e quasi completamente dimenticata. Si va avanti a furia di fortunate induzioni e di esperienze mentali.
 \end{quote}
 
Chiunque abbia dimestichezza con la filosofia di Enriques, riconosce la voce del maestro in questa citazione. Questi sconfinamenti filosofici di Severi non hanno serie conseguenze: Croce e Gentile sono troppo occupati a tenere d'occhio Enriques, che per loro, e non  a torto, costituisce la vera minaccia alla loro egemonia filosofica in Italia. 

Ma il  ``laissez-faire'' nei confronti di Severi si inceppa a causa di un suo nuovo tentativo di proporsi, o forse di imporsi (imitando e forse cercando di superare in questo come in altre cose Enriques), come filosofo della scienza, oltre che come matematico.  Il 24/11/1913 Severi pronuncia il discorso inaugurale \cite  {Sev1314} per l'anno accademico 1913--14 presso l'Universit\`a di Padova, intitolato ``Razionalismo e spiritualismo'' .  In esso, allineandosi del tutto ad Enriques, attacca Croce e Gentile (che a questo punto  erano gi\`a entrati in acerba polemica con Enriques) per la loro intransigenza filosofica, per la giustificazione storica e morale data da Croce alla Santa Inquisizione, per la loro sottovalutazione della scienza, per la loro intolleranza culturale. In particolare, Severi  si avventura a scrivere:

\begin{quote}
Io penso che la reazione contro il razionalismo e contro la scienza, in quanto ci allontana da siffatti principi umani, tende a riportarci verso abitudini mentali dannose. L'idealismo recente del Croce e del Gentile guida non alla correzione dell'errore ma alla 
persecuzione di chi erra.
\end{quote}

Benedetto Croce replica con il breve ma velenoso intervento \cite {Cro14} dal titolo \emph{E se parlassero di matematica?}, in cui, nell'accusare i ragionamenti di Severi di essere frutto di ``mentalit\`a  democratica (sic!) e massonica'',  scrive:

\begin{quote}
Al prof. Severi che \`e uomo di studio vorrei rivolgere una preghiera; ed \`e di non arrischiarsi a discutere concetti che appartengono a un campo a lui estraneo, e a entrare nel quale non so se abbia l'attitudine (ciascuno ha le sue attitudini), ma certo non ha la preparazione.
\end{quote}

L'ulteriore replica di Severi sul ``Giornale d'Italia'' del 5/2/1914 \`e anch'essa dura, sebbene espressa in tono ironico. Una posizione pi\`u equilibrata, meno settoriale e provinciale appare nel commento iniziale del redattore che indica come la contrapposizione tra Razionalismo e Spiritualismo sia da considerarsi superata:
\begin{quote}
A nostro umile parere, veramente di conflitto non si dovrebbe parlar pi\`u ...
\end{quote}
Un commento appropriato, che evidenzia come queste ``querelles'' siano decisamente di retroguardia rispetto allo sviluppo della filosofia scientifica a livello internazionale. Ad esempio nel 1908 Robert Musil  (1880--1942) discuteva presso l'Universit\`a di Berlino la sua tesi di dottorato in Filosofia dal titolo \emph{Sulle teorie di Mach} (cfr. \cite {Mus08}). Essa inizia cos\`i:
\begin{quote}
Ovunque, oggi, questioni metafisiche o gnoseologiche vengano sottoposte all'esame di una filosofia esatta, la parola dello scienziato ha un grande peso. Sono passati i tempi in cui l'immagine del mondo scaturiva in un atto di procreazione primordiale dal cervello del filosofo.
\end{quote}

Interessante comunque la difesa di Severi, nel citato articolo sul ``Giornale d'Italia'', dall'accusa da parte di Croce di essere massone:

\begin{quote}
... quanto al massonismo, il sen. Croce abbia la bont\`a di credere che non \`e il mio caso clinico! Se mi conoscesse di persona, egli non stenterebbe molto a convincersi ch'io non posseggo quel tanto di duttilit\`a morale e di diplomazia, che occorrono al perfetto massone.
\end{quote}
Se diamo fede, per\`o, all'annotazione da parte della polizia fascista al promemoria di Severi a Mussolini del 31/1/1929 (cfr. \cite[p. 104--107]{GuerNas05}), la ``duttilit\`a morale'' e la ``diplomazia che occorrono al perfetto massone'' se non le aveva ancora acquisite nel 1914, Severi le svilupp\`o  pi\`u tardi. O forse Croce, pur non conoscendo Severi ``di persona'', era ben informato.  

Fatto sta che, dopo questo episodio Severi cessa ogni forma di intervento filosofico diretto, e ogni polemica con Croce, di cui firmer\`a il manifesto degli intellettuali, e meno che mai con Gentile, con cui ebbe un rapporto di sodalizio umano e politico dagli anni '20 in  poi. E ci\`o va di pari passo con l'allontanamento personale da Enriques.

\section{Dopo la Grande Guerra}

Nel 1921, vincendo la concorrenza con Enriques, Severi viene chiamato all'Universit\`a di Roma sulla cattedra di Geometria Algebrica. Enriques verr\`a comunque anche lui  chiamato a Roma quasi contemporaneamente, grazie allo spostamento su altra cattedra di Guido Castelnuovo (1865--1952). A Roma Severi viene nominato rettore nel 1923, carica che lascia il 16 Novembre del 1925 con la convinzione, espressa nella sua lettera di dimissioni, ``[...] che ragioni politiche, menomando la fiducia che si riponeva in lui, hanno reso meno utile la sua opera di Rettore [...]''.  Le ``ragioni politiche'' da lui adombrate, sarebbero state la sua posizione avversa al governo a seguito del delitto Matteotti  nel 1924 e la sua adesione nel 1925 al ``Manifesto degli intellettuali antifascisti'' di Benedetto Croce. Tuttavia dalla corrispondenza tra Severi e Giovanni Gentile (cfr. \cite {GuerNas93}, si vedano in  particolare le lettere del 31/7, 8/8, 19/8 e 24/8/1925) si deduce che vi erano anche altre ragioni, di natura tecnico--amministrativa, per le sue dimissioni. Tra vari rilievi e attacchi subiti da Severi per la gestione dell'Ateneo romano, nel luglio del 1925 venne avviata un'inchiesta sulla situazione del bilancio del Policlinico. Severi attribu\`i attacchi e inchiesta proprio a ragioni politiche, scrivendo nella lettera del 31/7/1925:

\begin{quote}
Non escludo che possano dall'\emph{alto} esser venuti ordini di rendermi impossibile la vita rettoriale...  
\end{quote}

In una delle ultime lettere di Severi a Volterra del 20/8/1925, egli si riferisce all' ``affaire'' rettorato, e scrive:

\begin{quote}
E bisogna ch'io non mi discosti neppure di una piccola linea dalle disposizioni formali delle leggi e regolamenti, perch\'e si vigila di continuo su di me per cogliermi in fallo come avversario politico del ``regime''.
\end{quote} 

A partire dal 1925, con un'audace inversione di marcia, si avvicina al fascismo divenendone poi pieno sostenitore ed esponente tanto autorevole da potersi rivolgere direttamente a Mussolini. Fu lui, ad esempio, a suggerire nel 1929 a Mussolini, con un lungo promemoria (cfr. \cite[p. 104--107]{GuerNas05}), di imporre, come fu fatto nel 1931, il giuramento di fedelt\`a al regime ai professori universitari (a ci\`o si sottrassero solo una dozzina di docenti).  

Stranamente si iscrisse al Partito Nazionale Fascista soltanto nel 1932, rimanendo ad esso fedele fino alla sua caduta. In quegli anni Severi, profittando del suo  prestigio scientifico a livello internazionale, effettu\`o vari viaggi all'estero durante i quali, oltre a tenere conferenze e corsi di contenuto matematico,  si prodig\`o nel fare anche opera di propaganda a favore del regime. 

In  \cite {Sev59}, con una nota di tardivo pentimento per la sua compromissione col fascismo, riconobbe la sua scarsa accortezza in tema di politica,  affermando che:

\begin{quote}
... la matematica \`e l'arte di dare lo stesso nome a cose diverse: perci\`o i matematici sbagliano spesso quando si occupano di politica, giacch\'e la politica \`e invece l'arte di dare nomi diversi alle stesse cose.
\end{quote} 

Un'asserzione tipicamente Severiana questa, che  tende a fare del suo caso personale una teoria generale, da applicare a tutti i matematici. Per restare infatti ai matematici suoi contemporanei e a lui molto vicini, vista la colleganza che li univa presso l'Universit\`a di Roma, resta il fatto che Vito Volterra fu tra i pochi che non  giurarono fedelt\`a al fascismo cui fu sempre  avverso, e Tullio Levi Civita (1873--1941), simpatizzante delle idee socialiste e internazionaliste, rimase, al contrario di Severi,  sempre fermamente pacifista.  

Diviene Accademico d'Italia dal 1929 al 1944/1945\footnote {Dopo la caduta del regime fascista il 25 luglio 1943, l'Accademia d'Italia fu prima trasferita a Firenze, poi nella Villa Carlotta presso Tremezzo (lago di Como), dove, al tempo della Repubblica Sociale Italiana, sopravvisse fino al 25 aprile 1945, bench\'e soppressa dal governo legale fin dal 28 settembre 1944.}, vincendo di nuovo, sul filo del traguardo, la concorrenza di Enriques. La già citata probabile adesione alla massoneria, aggiunta alla circostanza che Enriques era ebreo, potrebbe aver giovato a Severi nella sua competizione per divenire accademico d'Italia. Va in ogni caso osservato che nessun ebreo divenne mai membro dell'Accademia d'Italia. 

Nel 1938, anno della promulgazione delle leggi razziali, fond\`o l'Istituto Nazionale di Alta Matematica di Roma,  che gli \`e attualmente intitolato.

Dopo la seconda guerra mondiale, Severi venne sottoposto a procedimenti di epurazione come docente universitario, come Presidente dell'Istituto di Alta Matematica e come socio della restaurata Accademia dei Lincei. Si difese affermando che, pur avendo aderito al fascismo, non aveva mai personalmente danneggiato nessuno. In tutti i procedimenti cui fu sottoposto venne prosciolto, mantenendo la presidenza dell'Istituto di Alta matematica e venendo rieletto Socio Linceo nel 1948 a seguito di un'amnistia generale.

Gli ultimi anni di vita di Severi, specie dopo la morte della moglie che per lui costitu\`i un duro colpo, vedono un ulteriore radicale cambiamento di rotta personale e politica, con un avvicinamento alla fede cattolica testimoniato dal gi\`a citato saggio autobiografico \cite {Sev59}.  

\subsection{La matematica}

I temi di ricerca di Severi non cambiano certo radicalmente dopo la Grande Guerra, anche se si possono riscontrare degli elementi innovativi, sia evolutivi che involutivi rispetto al periodo anteguerra. Questo periodo  \`e  caratterizzato da idee e proposte originalissime, ma spesso perseguite senza il dovuto rigore n\`e i necessari approfondimenti che dessero loro il corpo di compiute teorie scientifiche.  Infatti nel dopoguerra Severi, pur non perdendo la sua grande vena ispiratrice, la sua originalit\`a di vedute e producendo qua e l\`a delle vere perle, ha pi\`u \emph{congetturato} (per usare un eufemismo che nasconde varie affermazioni di dubbia esattezza) che \emph{dimostrato}.

Probabilmente i suoi lavori pi\`u rifiniti di questo periodo non sono in geometria algebrica ma in analisi complessa: si tratta di una serie di lavori dell'inizio degli anni '30 (cfr. \cite{Sev29, Sev31, Sev31/2, Sev31/3, Sev32, Sev32/2, Sev32/3}) dedicati ad un tentativo di sistemazione della teoria delle funzioni olomorfe di pi\`u variabili complesse, di cui Severi vedeva i profondi nessi con gli aspetti analitici e topologici della teoria delle superficie e variet\`a algebriche a pi\`u dimensioni. Per notizie pi\`u dettagliate in merito, si veda \cite [pp. 16--18]{Ves04}. 

I lavori di geometria algebrica riguardano per lo pi\`u innumerevoli, e sbagliati, tentativi di dimostrare  il teorema di completezza della serie caratteristica (cos\`i come del resto fece Enriques), estensioni del teorema della base a variet\`a di dimensione superiore, tentativi di ritornare sui punti ancora insoluti nella ``sistemazione'' da lui proposta della geometria numerativa. Questi lavori sono caratterizzati da una crescente ``fiducia in se stesso'', che spesso ha portato Severi a formulare asserzioni ingiustificate prendendole come teoremi, e fondandovi sopra teorie crollate quali castelli di carta di fronte alle obiezioni di critiche serrate.  In questo senso sono paradigmatici i molteplici lavori sulle serie di equivalenza, cui Severi si dedic\`o a partire dal 1932 per pi\`u di vent'anni. Si tratta di un argomento fondamentale (il concetto di \emph{equivalenza razionale} di cicli algebrici su una variet\`a) da lui per primo coraggiosamente e lungimirantemente intrapreso, ma nel quale, accanto a idee originali e a intuizioni geniali, produsse anche molti gravi errori ed incongruenze (cfr. \cite [\S 3] {BCP04}). 

David Mumford, in un'appendice alla ristampa del libro \cite {Zar71} di Oscar Zariski  (1899--1986) dice:

\begin{quote}
In the period 1935 to 1950, SEVERI published many papers on series of equivalence
and its generalizations to higher dimensions. It is hard to untangle everywhere
what he conjectured and what he proved and, unfortunately, some of his conclusions
are incorrect.
\end{quote}

Quel che \`e certo \`e che Severi non riusc\`i nemmeno a dare una definzione
precisa di {equivalenza razionale} che sarebbe dovuta essere alla base della sua teoria, al punto che nel Congresso Internazionale dei Matematici 
di Amsterdam del 1954 ``the
definition itself of series of equivalence was debated sharply'', come dice Mumford, l.c.: in realt\`a, a seguito delle obiezioni e degli attacchi di Pierre Samuel  (1921--2009) al cospetto di A. Weil, Severi non  riusc\`i neanche a terminare il suo intervento. 

Un giudizio meno tagliente \`e dovuto a William Fulton, nel suo libro \cite  {Ful84}:

\begin{quote}
It would be unfortunate if Severi's pioneering works in this area were forgotten;
and if incompleteness or the presence of errors are grounds for ignoring Severi's
work, few of the subsequent papers on rational equivalence would survive.
\end{quote}

In ogni caso, Mumford \`e stato sempre alquanto duro con Severi, forse anche a causa dell'avversione nutrita verso Severi dal suo maestro Zariski, molto vicino ad Enriques che invece con Severi ebbe un pessimo rapporto. Ad esempio,  nel recente articolo \cite {Mum11}, a commento dei tentativi di Enriques di ``aggiustare'' la sua dimostrazione del teorema di completezza della serie caratteristica (e alla sua polemica con Severi su questo argomento su cui torneremo), Mumford dice:

\begin{quote}
In short, Enriques was a visionary. And remarkably his intuitions never
seemed to fail him (unlike those of Severi whose extrapolations of known
theories were sometimes quite wrong).
\end{quote}

In altra occasione, commentando un fondamentale risultato da lui dimostrato (cfr. \cite {Mum69}), su cui invece Severi aveva a suo tempo asserito qualcosa di completamente falso, pur avendo delle tecniche che lo stesso Mumford usa per abbordarlo nel modo giusto,  dice, con non velata ironia:

\begin{quote}
Now, after criticizing Severi like this, I have to admit the following: the method
of disproof of (1)-(3) is due entirely to Severi: Severi created, in fact, a very excellent
tool for analyzing the influence of regular 2-forms on $F$ on his systems of
equivalence. One must admit that in this case the technique of the Italians was
superior to their vaunted intuition.
\end{quote}

\subsection{La politica e il versante accademico}

L'adesione di Severi al fascismo toccò un punto di non ritorno col citato promemoria a Mussolini del 1929, in cui egli formula il suggerimento di imporre ai professori universitari il giuramento di fedelt\`a al regime. Come non \`e raro in lui, egli presenta la cosa come dettata da un'esigenza generale di normalizzazione, ma in buona misura tale esigenza \`e ritagliata su se stesso: 

\begin{quote} ... poich\'e i nostri sentimenti di devozione alla Patria non posson non esser quei medesimi che nel 1914--15 fecero delle Universit\`a le fucine dell'intervento, i Professori non desiderano oggi che di poter cooperare con lealt\`a ed in un'atmosfera di fiducia, alla grande opera di ricostruzione, di cui il Capo ha posto i saldi piloni. [...]
Il mio caso \`e quello di tanti. Uscito nei primi mesi del 1915 dalle file socialiste (con una dichiarazione che, se ben ricordo, fu, almeno in parte, pubblicata nel ``Popolo d'Italia'') mi arruolai volontario allo scoppiar della guerra e fui sempre combattente al fronte. Finita la guerra, a Padova, dove ero Direttore di quella Scuola di Ingegneri, fronteggiai coi combattenti il movimento bolscevico.
\end{quote} 

Severi rimase fedele al regime anche dopo la promulgazione delle infami leggi razziali. E Severi \`e uno dei membri della Commissione Scientifica dell'Unione Matematica Italiana che delibera, il 10 Dicembre 1938\footnote {Gli altri membri presenti erano Luigi Berzolari (1863--1949), Enrico Bompiani (1889,--1975) , Ettore Bortolotti (1866--1947), Oscar Chisini (1889--1967), Annibale Comessatti, Luigi Fantappi\'e (1901--1956), Mauro Picone (1885--1977), Giovanni Sansone (1888--1979), Gaetano Scorza (1876--1939), assente giustificato Leonida Tonelli (1885--1946).}, l'approvazione di un comunicato  che afferma, tra l'altro, che ``la scuola matematica ita\-liana, che ha acquistato vasta rinomanza in tutto il mondo scientifico, \`e quasi totalmente creazione di scienziati di razza italica (ariana)''.  Si tratta di un atto di grave compromissione col regime fascista  di cui molto probabilmente Severi fu uno degli ispiratori, un atto tanto pi\`u grave quanto rivolto  al meschino obiettivo di far s\`i che ``nessuna delle cattedre di matematica rimaste vacanti in seguito ai provvedimenti per l'integrit\`a della razza, venga sottratta alle discipline matematiche''. Nonostante ci\`o, e la posizione di indiscusso predominio nella matematica italiana che l'espulsione degli ebrei consegn\`o definitivamente a Severi, questi riusc\`i, nel secondo dopoguerra a ``farsi perdonare'' la sua totale adesione al regime. Ad esempio, leggiamo in \cite {GooBab11}:

\begin {quote}
Among the first, if not the first, to come to Severi's defense in the post-World War II era was Beniamino Segre, who had lost his own academic chair in Bologna in 1938 in the wake of the government's anti--Jewish legislation. Relieved, at Severi's behest, of his duties as an editor of Italy's oldest scientific journal in the wake of the regime's racial laws, Segre nevertheless insisted that rumors 
of Severi's ``supposed anti-Semitism'' were ``flimsy'' and based on ``some misunderstanding''.
\end{quote}

In effetti l'attaccamento di Segre al Maestro ha del singolare, e  perfino dei tratti patetici. Basti tenere presente la lettera del 2/1/1932, di Severi a Segre, riportata in \cite{GooBab11}, in cui il primo, nel criticare aspramente uno scritto del secondo, fa una smaccata apologia di se stesso, al contempo impietosamente ed ingiustamente sminuendo il lavoro di Corrado Segre, suo antico maestro, nonch\'e quello di Castelnuovo ed Enriques, che \`e ridicolmente ridimensionato. Citiamo per esteso dalla lettera, nella traduzione di \cite{GooBab11}. Citazione che ci pare di particolare interesse anche per delucidare il carattere di Severi e la visione che egli aveva del suo lavoro in relazione a quello degli altri, che si manifestano con brutale franchezza dato la forma privata della corrispondenza:

\begin {quote}
My dear Segre,

The general outline [of B. Segre's draft] is mediocre in several places,
especially where you talk about algebraic geometry.
It lacks perspective so that a reader who doesn't know much will not be able to understand the hierarchy of ideas
and names.\\
1) The work of C. Segre has been overrated ...
Segre, for example, did not
prove anything major in the field of
geometry of curves although he did
carry out a very significant revision of
the subject. His contributions to higher
dimensional projective geometry are
overrated when compared, for example,
with those of Veronese. This exaggerated
evaluation is probably explained
by your love of him as a disciple.\\
2) The work by Veronese is underrated.
In Italy he was the true creator of higher
dimensional projective geometry.\\
3) The work of Castelnuovo has been
overrated as has been that of Enriques,
especially when compared to that of
[Max] Noether, whose name you have
completely neglected in your discussion
of surfaces.\\
4) My work has been underrated, which
seems odd to me since you were my
student, and, in addition, your affection
for your first teacher [Segre] has caused
you to overestimate his work.
[...]
Beginning in 1904 I developed new
ideas that untied the Gordian knots
that had been bound so tightly up to
that time. I myself untied most of them,
such as the characterization of irregular
surfaces from both the transcendental
and algebraic point of view. Even
setting aside my work on conceptual
clarifications as well as my work on the
hyperelliptic surfaces with Enriques,
which I am willing to do, this does not
justify the humiliating description of
my status as arrived' that you have
written on page 12, thus putting my
work and that of Castelnuovo-Segre's
[Severi may have meant Enriques here
ed.] on two completely different levels.
You need to weigh your words!

Also forgotten by you were my Theorem
of the Base [Sev06] and my work
on the geometry of varieties of higher
dimension. How could you have done
this when you discuss Italian geometry?
In addition, there is no mention of the
fact that I am the only one among living
Italian algebraic geometers who has
created a school.

You have also underestimated my contributions
to enumerative geometry. If
only you could understand them. All
of this is not to reproach you, because
you certainly have done your best. Although
your mathematical knowledge
is wide, you currently do not have a
deep enough understanding in the vast
field of algebraic geometry to allow you
to have a reliable perspective on the
subject. But I am also surprised that
the comparative evaluations that we
discussed many times in the past did
not have an effect on you even though
I was always very conscientious about
being objective.

\end{quote}

Ci\`o nonostante,  nel 1962, nella sua commemorazione di Severi \cite {Seg62}, Segre afferma, con una buona dose di esagerazione:
\begin{quote}
FRANCESCO SEVERI fu un gigante del pensiero matematico moderno. La sua opera scientifica, geniale e profonda, sotto vari aspetti decisiva per i progressi recenti e futuri dell'Algebra e della Geometria, \`e ricchissima d'idee luminose e feconde, oltrech\'e di risultati fondamentali e stupefacente anche per ampiezza e multiformit\`a. [...]
La proteiforme e maestosa figura di FRANCESCO SEVERI ha dominato con le sue doti affascinanti il mondo scientifico e culturale italiano per oltre mezzo secolo.
\end{quote} 

Affermazioni che peraltro evidenziano un certo provinciale distacco dai grandi sviluppi che la matematica, e in particolare la geometria algebrica, cui si assisteva in quegli anni ad opera delle scuole francese, americana e russa.

Tornando alle leggi razziali, e volendo delucidare ulteriormente il carattere di Severi e il suo atteggiamento in quella triste occasione,
non si pu\`o dimenticare la lettera di Severi a Tullio Levi Civita (1873--1941) del 18/10/1938\footnote{Si conservano poche lettere di Severi a Levi Civita, quattro del 1921, tre del 1929, una del 1938 (con risposta  da parte di Levi Civita), una del 1940. Ringraziamo Pietro Nastasi per averci fornito tali lettere.}:

\begin{quote}
 Carissimo Tullio,
In obbedienza alle superiori disposizioni trasmessemi dal Ministero dell'Educazione nei riguardi degli ``Annali di Matematica'', debbo procedere, di concerto con la Casa Zanichelli, alla tua sostituzione della Direzione del periodico.
Prima di farlo, tanto io che il Dott. Della Monica, abbiamo desiderato di dartene notizia e di ringraziarti vivamente della cooperazione efficace ed altamente autorevole che in tanti anni hai dato come Condirettore degli ``Annali''.
Questa separazione non muta in nulla i nostri legami di buona amicizia personale, che dura inalterata da pi\`u di un trentennio. [...]
 \end{quote}
 
 Come si vede non una parola di rincrescimento n\`e di solidariet\`a.  La serena lettera di risposta di Levi Civita non sembra mostrare elementi di astio.

 \begin{quote}
Carissimo Francesco,
Ricevo a Roma, dove sono rientrato da poco pi\`u di una settimana, la tua lettera del 18, diretta a Padova. Nel prenderne atto, ricambio con pari amicizia gli affettuosi sentimenti che hai voluto esprimermi in modo lusinghiero e cordiale.
Una stretta di mano.
 \end{quote}

Severi e Levi Civita erano legati  da antica amicizia e colleganza maturata nel periodo padovano: come vedremo, Severi stesso riconosceva, Levi Civita si era battuto per la sua chiamata a Padova, cos\`i come si batt\`e per la sua chiamata a Roma, e si rallegregr\`o per la sua nomina ad Accademico d'Italia. Un'amicizia cementata, nell'epoca patavina, anche da idee politiche collimanti, idee cui Levi Civita rimase fedele per tutto il corso della sua vita, mentre, come abbiamo visto, Severi rinneg\`o. 

Con la chiamata a Roma, l'adesione al fascismo, la nomina ad Accademico d'Italia, Severi va assumendo nel corso degli ani '20 del XX secolo sempre pi\`u un ruolo di preminenza nella matematica italiana. Ruolo che prima di lui era stato ricoperto da Volterra, il quale, a causa delle sue posizioni politiche apertamente liberali e sempre ostili al fascismo (anche lui, come Severi, fu tra i firmatari del manifesto Croce, ma non rinneg\`o mai tale scelta), cadde lentamente in disgrazia, perdendo nel 1927 la presidenza del CNR da lui creato nel 1923, poi venendo allontanato dall'Universit\`a nel 1931 e succesivamente dall'Accademia dei Lincei, di cui era stato presidente dal 1923 al 1926, causa il rifiuto di giurare fedelt\`a al regime. 

Il ruolo di preminenza di Volterra si basava non solo sulla assoluta eccellenza della sua produzione matematica teorica, che ne fa uno dei fondatori della moderna analisi funzionale, ma anche sulle sue doti di  matematico applicato, uno scienziato moderno e a tutto tondo, in grado di risolvere problemi della pi\`u varia natura, adattando all'uopo con grande flessibilit\`a raffinati strumenti matematici. Tra questi, basti ricordare il modello  \emph{preda--predatore}, oggi noto come \emph{modello Lotka--Volterra}\footnote{Alfred James Lotka (1880--1949) fu un matematico, statistico e chimico fisico, si occup\`o di dinamica delle polpolazioni, proponendo il modello in questione contemporaneamente a, ma indipendentemente da, Volterra.}. 

La personalit\`a  poliedrica di Volterra era tale da porre in sordina l'eterna, e purtroppo sempre attuale, diatriba di retroguardia  tra ``matematica pura'' e ``matematica applicata''. Una diatriba che, nel passato come nel presente, si fonda per lo pi\`u sul tentativo di affermazione dell'una sull'altra, come fossero discipline distinte o peggio in conflitto, ed ha come unico vero obiettivo la mera conquista di fette di potere e di finanziamenti. 

Tramontato l'astro di Volterra, la diatriba scoppi\`o vivace alla fine degli anni '30, e vide coinvolti da un lato Severi, assertore della preminenza della ``scienza pura'', dall'altro Mauro Picone, fascista della prima ora e fondatore dell'Istituto Nazionale delle Applicazioni del Calcolo del CNR, propugnatore, accanto al  valore teorico della matematica,  anche dell'importanza culturale e sociale delle sue applicazioni. Del conflitto, tutto interno alla  politica culturale e accademica fascista (entrambe le parti cercavano di piegare le roboanti asserzioni generali della propaganda di regime ai loro punti di vista), si espresse con toni civili nella forma, ma che non nascondevano una decisa lotta di potere accademico. Di essa offre  un ottimo sunto Pietro Nastasi in \cite[\S 6.3] {Nas06} e ne d\`a conto anche un osservatore del tempo con nessuna simpatia per il regime,  Francesco Tricomi (1897--1978), in un brano riportato in \cite [pp. 81--82] {GuerNas93}:

\begin{quote} \`E curioso osservare come il Consiglio delle Ricerche, che cominci\`o a funzionare effettivamente intorno al 1935, fino alla guerra ebbe s\`i un Comitato per la matematica applicata, ma non uno per la matematica senz'altro\footnote{In effetti da Comitato per la matematica del CNR fungeva l'Unione Matematica Italiana, creata nel 1922 per iniziativa di Volterra, cfr. \cite {Cil17}.}. La ragione (naturalmente, non detta) fu di evitare che se ne impadronisse Severi che, divenuto inopinatamente Accademico d'Italia nel 1929 e altrettanto esuberante fascista quanto prima era stato antifascista, voleva essere (e, in certa misura, fu) il ``padrone'' della matematica italiana nel periodo fascista. 
\end{quote}

Accanto a questo conflitto  si colloca il contrasto con Enriques, cui pi\`u volte abbiamo fatto cenno. Severi non era nuovo a polemiche che, sorte in ambito scientifico,  scadevano poi in conflitti apertamente personali. Gi\`a prima della Grande Guerra, ricordiamo quella con Michele De Franchis (cfr. \cite [pp. 21--22] {CilSer91}) e quella esplosa proprio durante la guerra, con Giovanni Zeno Giambelli (1879--1935) su questioni metodologiche riguardanti l'approccio alla geometria numerativa (cfr. \cite [\S 2.6.1] {BriCil88} e \cite [S 2.6] {BriCil95}). Nonostante di Giambelli venga oggi internazionalmente riconosciuto il valore per il suo approccio aritmetico alla teoria della intersezione (cfr. \cite {Lak94}), egli, anche a causa del suo contrasto con Severi, non riusc\`i mai ad ottenere un posto di professore universitario. 

Il contrasto con Enriques si pone anch'esso sulla linea di confine tra l'ambito scientifico, quello della politica accademica e quello personale: si tratta in verit\`a pi\`u di un conflitto per una supremazia culturale  che di un vero e proprio dissidio su argomenti scientifici i cui prodromi, come abbiamo visto, si manifestano gi\`a prima della Grande Guerra. 

Passando al dopoguerra e tralasciando l'episodio della chiamata a Roma, un momento chiave dell'approfondirsi di un solco tra i due, \`e da collocarsi nell'episodio del 1921 della segnalazione da parte di Severi in \cite [Articolo 5, nota (15)] {Sev21/2} di un grave errore nella dimostrazione di Enriques del teorema di completezza della serie caratteristica del 1904, che nel frattempo era stato mirabilmente dimostrato per via analitica da Henri Poincar\'e (1854--1912) in  \cite {Poi10}, dimostrazione che viene in parte semplificata da Severi in \cite {Sev21}. In verit\`a Severi garbatamente accenna:
\begin{quote}
[...] ad un punto del procedimento geometrico con cui si dimostra la completezza della serie caratteristica che abbisogna di qualche ulteriore indagine, la quale potr\`a eventualmente portare a limitazioni inessenziali, lasciando per\`o integra la sostanza del fatto.
\end{quote}

La reazione di Enriques a questa cauta critica giunse  pi\`u tardi in una serie di note del 1936--37. Severi rispose con due memorie del 1941--42, e B. Segre si inser\`i con due lavori del 1938--39 (la vicenda \`e narrata, dal punto di vista di Enriques, nel suo libro postumo \cite  [Cap. IX, \S 6]  {Enr49}). Con il passare degli anni, la polemica assunse toni assai aspri, scadendo oltre il limite del lecito, senza  peraltro concretizzarsi in un vero dissidio scientificamente significativo: Enriques, come rilevato da Mumford, forse aveva le idee giuste, ma non i mezzi tecnici per realizzarle, e Severi, nelle sue risposte, aggiungeva errori ad errori.

Ma nel frattempo, alla fine degli anni '20 tra i due si manifestavano anche su altri fronti dissidi, conflitti e veri e propri litigi personali, tali da rendere i loro rapporti  pi\'u che difficili, impossibili. Cause degli attriti: la nomina a socio dell'Accademia d'Italia, di cui si \`e gi\`a detto, la collocazione editoriale e l'adozione dei rispettivi testi scolastici (cfr. \cite{GiacTea15, Lin04}), la mancata collaborazione di Severi all'edizione della Enciclopedia Italiana, alla direzione della cui sezione matematica Gentile aveva chiamato Enriques, che intanto aveva con lui ricucito i rapporti dopo gli scontri filosofici  dell'anteguerra \cite {Bol04}. Ad esempio \`e del 24/5/1928 una lettera di Severi a Gentile in cui si dice:

\begin{quote}
[...] con un uomo come Enriques [...] che io giudico il pi\`u inadatto per compiere un'opera di valutazione obiettiva, come quella che si richiede in un'Enciclopedia, io non posso avere pi\`u nulla in comune e tanto meno relazioni di quasi subordinazione.
\end{quote}

Va ricordata infine l'aspra polemica che segue la recensione \cite {Enr34} di Enriques delle \emph{Lezioni di Analisi} \cite {Sev33} di Severi, in cui Enriques, che evidentemente non si asteneva dal cercare terreni di ulteriore scontro, attacca Severi per il suo eccessivo rigorismo nella trattazione del cosiddetto \emph {Criterio di Pl\"ucker--Clebsch}. Essa dar\`a la stura ad una serie di articoli, dell'uno e dell'altro, che difendendo le proprie posizioni, attaccano l'avversario. 

Infine \`e un fatto che Severi, a partire dagli anni Venti, quando appunto avvenne il trasferimento a Roma, cerc\`o di assumere un ruolo guida sia accademico, sia culturale, innanzitutto nella geometria algebrica, poi pi\`u in generale nella matematica nazionale. Ci\`o si svilupp\`o in varie direzioni, e una di queste fu l'iniziativa di provare a riscrivere la geometria algebrica italiana dal suo punto di vista. In questa luce va vista la pubblicazione delle \emph{Vorlesungen} \cite {Sev21} del 1921 e pi\`u ancora  il tentativo di scrittura del \emph{Trattato di Geometria Algebrica} \cite {Sev26}, un progetto molto ambizioso (troppo, perfino per la grande determinazione e le capacit\`a lavorative di Severi) che avrebbe dovuto costituire una vera e propria enciclopedia dettagliata della materia, ma di cui non vide la luce che il solo primo volume nel 1926. Come gi\`a osservato in  \cite [\S 2.6] {BriCil88}, questo progetto non poteva certo riuscire gradito a Enriques, che un 
simile monumentale progetto aveva concepito e in gran parte realizzato con le \emph{Lezioni sulla teoria geometrica delle equazioni e delle funzioni algebriche} scritte in collaborazione con l'allievo O. Chisini \cite {EnrChi}. E infatti Severi, nel difendersi dagli attacchi di Enriques alle sue \emph{Lezioni di Analisi}, attacc\`o a sua volta in \cite {Sev34}  le \emph{Lezioni} per lo scarso rigore nella trattazione proprio del  \emph {Criterio di Pl\"ucker--Clebsch}. Le riserve di Enriques, oltre che dettate dalla difesa di una posizione di preminenza scientifica, avevano anche una motivazione metodologica, come si rileva dalle critiche da lui mosse alle \emph{Lezioni di Analisi}. \`E forse questo l'unico punto in cui ci sia del succo scientifico nel contrasto tra i due. Enriques postulava infatti l'intuizione come prevalente metodo di azione scientifica ed era disinteressato, se non contrario, ad una esposizione sistematica dei concetti. Severi provava invece, forse anche sulla spinta delle prime critiche metodologiche provenienti da Oscar Zariski (1899-1986), da A. Weil e dall'incipiente Bourbakismo, a costruire un sistema completo e rigoroso (quale solo pi\`u tardi riusc\`i, e con ben altri strumenti, a Jean Paul Serre, Alexander Grothendieck (1928--2014), ed altri). Fu ben lontano dal riuscirci per mancanza, ed anche per rifiuto, di opportuni strumenti algebrici, topologici, ecc. Sfortunatamente si tratta di un dissidio di retroguardia, fin troppo personalizzato, tutto interno ad una scuola ormai in una decadente fase di ripiegamento su se stessa.
 
 \section{Conclusioni}  
 
 Severi \`e un personaggio complesso e problematico, a tinte taglienti, in cui luci ed ombre convivono in modo inestricabile, dotato di prorompente ambizione, di grandissimo ingegno e forza di lavoro, ma altrettanto privo di equilibrio di giudizio e autocritica, offuscato com'era da un invincibile, autoreferenziale egocentrismo. Tendenza che si amplifica col passare degli anni, presentando dei punti critici in occasione di particolari successi della sua carriera, come la chiamata a Roma o la nomina ad accademico d'Italia. La Grande Guerra non pare di per s\`e aver influito radicalmente n\`e sulla sua visione della matematica n\`e, in modo decisivo, sui suoi interessi all'interno di questa. Tuttavia essa, o gli anni che immediatamente la seguono, costituisce un giro di boa della sua vita scientifica, accademica e personale. Con la Grande Guerra si esaurisce il periodo pi\`u fertile dei suoi risultati, per dar posto pi\`u a brillanti intuizioni, accompagnate da sconcertanti errori e cadute di stile, che a teoremi. Anche in questo ambito si vede il progressivo manifestarsi della sua tendenza all'egocentrismo e alla autoreferenzialità che lo portano a nutrire una eccessiva fiducia nella validità delle sue intuizioni. Al volgere della Grande Guerra, la sua visione politica cambia radicalmente, forse anche per motivi opportunistici, sebbene il personaggio sia tanto complesso che attribuire solo a questo i suoi cambi di campo sia forse un po' riduttivo. Certamente il suo coinvolgimento col regime fascista si spiega con la sua grande ambizione di ottenere, mediante l'ossequio al potere, anche i riconoscimenti che reputava meritasse. Ma, come abbiamo cercato di porre in evidenza, c'\`e di pi\`u, e precisamente anche l'idea di poter in qualche modo piegare lo stesso potere fascista alle sue idee ed esigenze, che nel suo egocentrismo egli tendeva a  reputare di indole generale. Esempio tipico, il riuscire a lavare, mediante il giuramento di fedelt\`a al regime, il peccato originale dell'antifascismo. C'era poi il concetto di portare, facendo pesare il suo personale prestigio e la sua influenza, acqua al mulino della ``scienza pura'', della matematica teorica in particolare, difendendola, anche da un punto di vista epistemologico,  dalle derive utilitariste e applicative.  Per contro, con la Grande Guerra il suo incipiente cammino filosofico si interrompe bruscamente. \`E ben chiaro che lo scontro con i filosofi idealisti Croce e Gentile, che sempre pi\`u, anche se da sponde politiche opposte, stanno prendendo il predominio culturale in Italia, non gli sembri la cosa pi\`u opportuna in vista della agognata sua scalata alle vette del potere accademico.  
  
Ulteriori auspicabili ricerche (ad esempio, non tutti i contributi matematici di Severi sono stati analizzati nel giusto dettaglio, della sua corrispondenza privata ci \`e giunto ben poco, ulteriori indagini negli archivi universitari, ministeriali, militari, di partito, sono forse opportune) potranno ulteriormente illuminare questo notevole personaggio. 

\section{Appendice}\label{App}
Lettera di Severi al giornale ``L'Adriatico'' del 9/3/1915: 

\begin{quote}
Onorevole sig. Direttore, accolgo volentieri il suo cortese invito a partecipare alla discussione che l'``Adriatico'' ha aperto sui doveri degli italiani, e specialmente dei veneti, in quest'ora tragica e decisiva della storia, anche perch\'e a me, uomo di studio pi\`u che uomo di parte, sar\`a forse possibile di scriver parole che possan essere valutate serenamente da' miei compagni del partito socialista ufficiale, come dagli avversari, e che pertanto contribuiscano, sia pure in modestissima misura, a ricondurre a quel rispetto reciproco delle idee e a quella concordia degli animi, di cui qui nel Veneto, ancor pi\`u che altrove, c'\`e tanto bisogno.
Io spero e credo che l'atteggiamento degli organi direttivi del mio partito, in questo grave momento, sia l'espressione dello stato di angoscioso dissidio in cui ogni socialista d'intelletto e di cuore si trova fra gli imperativi ideali della propria fede e la percezione delle necessit\`a ineluttabili dell'ora presente; piuttosto che indice di un proposito d'azione maturato e metodicamente perseguito.
Ma se cos\`i \`e, e se \`e pur vero, secondo io penso, che il partito socialista, come organismo politico, non potrebbe mai farsi promotore di un intervento guerresco, assai meglio parmi si provvederebbe, se la protesta socialista contro la guerra fosse, in ogni occasione, contenuta nel campo puramente ideale, riconoscendo nello stesso tempo la ineluttabilit\`a d'una situazione che non ci è dato oggi di modificare, appunto perch\'e deriva da condizioni sociali che il partito nostro non pu\`o cambiare di colpo.
Porsi da un punto di vista di assoluta negazione di problemi che esistono e che reclamano una soluzione indifferibile, significa lasciarsi cullare dalla ingenua illusione di poter violentare lo svolgersi dei fenomeni storici, e venir quindi, in ultima analisi, a contraddire a quello che \`e lo spirito animatore della dottrina socialista.
Un atteggiamento meno assoluto della Direzione del nostro Partito, sarebbe importantissimo anche dal punto di vista politico, giacch\'e lascerebbe ad ogni inscritto la libert\`a di valutare gli elementi reali della situazione, secondo la propria coscienza di cittadino italiano, e nello stesso tempo consentirebbe ad ognuno di noi di continuare la propaganda socialista fra le masse, additando loro quali disastri immani conseguono dall'ordinamento capitalistico della societ\`a.
Io, che sono convinto della necessit\`a dell'intervento dell'Italia a fianco della Triplice Intesa, sento di non aver mai provato un odio cos\`i implacabile contro la guerra -- la quale non crea, ma sfrutta valori morali gi\`a esistenti; -- n\'e di aver mai desiderato, con altrettanto ardore, profondi rinnovamenti sociali, come da quando assistiamo alla spaventosa ecatombe di vite umane, all'enorme distruzione di ricchezza, all'acutizzazione del disagio economico del proletariato, al dispregio del diritto e delle bellezze dell'arte che la guerra europea trascina con s\'e.
Ma come non v'ha uomo cui la violenza ripugni, che ad essa non possa contro ogni sua voglia essere costretto; come non v'ha socialista che, vivendo e vestendo panni in questa societ\`a borghese, non s'adatti, nelle pratiche contingenze della vita, a ci\`o che l'ambiente gli impone senza che per questo egli rinunci a dar l'opera sua per un migliore domani, cos\`i non trovo vi possa essere contraddizione sostanziale fra la fede nei nostri ideali e l'azione che oggi cagioni storiche superiori alla nostra volont\`a possono prescriverci.
Vi sar\`a \`e vero, per chi ami dilettarsi in cos\`i tragico momento di quisquilie dialettiche, una contraddizione formale; ma sciaguratamente le pi\`u angosciose situazioni sentimentali si sciolgono di rado alla stregua della logica pura.
Eppoi il partito socialista non ha forse riconosciuto che nella pratica quotidiana conviene adattarsi ad un programma minimo e non evitare talvolta contatti con le frazioni pi\`u illuminate della borghesia, quando occorra, ad esempio, contrastare la vittoria di partiti i quali minaccino di prevalere in modo pericoloso per le libert\`a politiche, che costituiscono il presupposto delle conquiste economiche del proletariato? E perch\'e dovremmo racchiuderci in una formola d'intransigente negazione, proprio in una questione che di gran lunga trascende la importanza della minuscola politica d'ogni giorno, e che \`e in fondo ancora una questione vitale di libert\`a?
Giacch\'e \`e ben vero che le cause di questa guerra sono giustamente capitalistiche, ma non si pu\`o disconoscere che, sia per le brutali violazioni del diritto naturale dei popoli compiute dalla Germania, sia per l'esistenza di molte questioni insolute, sia infine per l'interesse di alcuni belligeranti, e soprattutto dell'Inghilterra, affinch\'e vengano rispettate le nazionalit\`a minori (``L'interesse e il dovere spingono l'Inghilterra nella stessa direzione'', hanno scritto i professori dell'Universit\`a di Oxford), la guerra \`e andata acquistando, in modo prevalente, il carattere d'un conflitto fra due opposte concezioni dei diritti e delle forze, che debbono prevalere nel mondo moderno. Inoltre, secondo la lettera e lo spirito della dottrina marxista, il socialismo potr\`a e dovr\`a succedere agli attuali ordinamenti, soltanto allora che la civilt\`a sia passata per tutte le fasi del suo sviluppo, tra le quali vi \`e appunto la conquista delle unit\`a e delle autonomie nazionali. Di guisa che, per dirla con una frase scritta in questi giorni nell'``Avanti'' da Enrico Leone, la Nazione diventa la porta d'ingresso dell'Internazionale.
E quando si parla della Nazione non ci si appiglia ad un ``diversivo borghese'', poich\'e la Nazione \`e una formazione storica naturale, la quale vive nelle tradizioni di lingua, di arte, di cultura, di ciascuna razza, e sta al disopra e al di fuori delle iniquit\`a degli ordinamenti statali.
Riconosco come vi siano purtroppo molti, i quali, per le condizioni di inferiorit\`a culturale e materiale in cui si trovano, non certo per colpa loro, non possono sentire tutto il valore spirituale dell'idea di Nazione; ma essi comunque non dovranno disconoscere che il dominio straniero rappresenta sempre un altro sfruttamento, da Nazione a Nazione, che viene ad aggiungersi allo sfruttamento del capitalista sul salariato.
Eppoi in qual modo si concreterebbe la solidariet\`a internazionale se, fino a quando non sar\`a pi\`u diffusa la coscienza della disastrosa follia degli armamenti, di fronte a tentativi di sopraffazione imperialistica a danno di altri popoli, non si fosse disposti anche a sacrifici di sangue?
D'altronde i socialisti, predicando l'avversione alle spese militari, hanno sempre presupposto la sincerit\`a e l'efficacia della propaganda antimilitarista negli altri paesi, ed hanno inteso con ci\`o di cercar di diminuire la possibilit\`a di conflitti armati fra i popoli, ma non gi\`a di negare le idealità nazionali. Allorch\'e la Patria sia in pericolo, ancor pi\`u impellente sorge quindi per noi socialisti il dovere di difenderla, avvalorando agli occhi di chi ci considera utopisti, la nostra persuasione che dalla coscienza di un buon diritto possa -- ove occorra -- sprigionarsi la pi\`u grande delle forze.
Ed io credo per certo che sul riconoscimento di questo dovere, la stragrande maggioranza dei socialisti italiani sia senza esitanze concorde, anche se qualche eccesso polemico possa a taluno far supporre il contrario. Non imprigioniamoci dunque nell'adorazione di formule assolute, giacch\'e il pericolo per il nostro Paese \`e insito nella grave situazione internazionale, la quale potrebbe trascinarci pi\`u tardi, anche nolenti, ad una guerra disastrosa per l'Italia e pi\`u particolarmente per il Veneto.\\

Prof. Francesco Severi, dell'Universit\`a di Padova.\\

P.S. -- Ragioni varie hanno fatto ritardare per circa due settimane la pubblicazione di questa mia lettera. Non ho ora nulla da mutare, ma di fronte al fatto -- segnalato anche ieri in queste colonne dall'amico Gino Piva -- che le condizioni economiche del proletariato veneto vanno di giorno in giorno aggravandosi, in modo veramente doloroso e allarmante, desidero di aggiungere una parola di viva deplorazione per l'inerzia del Governo, il quale sembra non abbia capito e non capisca che la preparazione non deve limitarsi alle sole provvidenze militari. Come si potrebbe sperare che le masse popolari offrissero la necessaria resistenza morale e materiale, se la nostra regione dovesse essere assoggettata, dall'intervento dell'Italia nel conflitto europeo, ad altre e ben più dure prove?
Provvedimenti eccezionali (lavori e sovvenzioni dello Stato ai Comuni) urgono qui nel Veneto per fronteggiare la grave crisi. Altro che proibire i comizi!
\end{quote}

\end{document}